\documentclass{amsart}
\usepackage{amssymb,latexsym,epsfig,color}

\newtheorem{theorem}{Theorem}[section]
\newtheorem{proposition}[theorem]{Proposition}
\newtheorem{corollary}[theorem]{Corollary}
\newtheorem{conjecture}[theorem]{Conjecture}
\newtheorem{lemma}[theorem]{Lemma}

\newenvironment{numtheorem}[1]{\medskip\noindent {\bf Theorem~ #1.}\em}{}

\theoremstyle{definition}
\newtheorem{definition}[theorem]{Definition}
\newtheorem{remark}[theorem]{Remark}
\newtheorem{example}[theorem]{Example}

\def\proof{\smallskip\noindent {\it Proof: \ }}
\def\endproof{\hfill$\square$\medskip}
\def\Z{\mathbb{Z}}
\def\N{\mathbb{N}}
\def\m{{\bf{m}}}
\def\Sfield{S}

\newcommand{\Ass}{\mbox{\upshape Ass}\,}

\newcommand{\Gin}{\mbox{\upshape Gin}\,}
\newcommand{\field}{{\bf k}}
\newcommand{\Field}{{\bf  K}}

\newcommand{\mult}{\mbox{\upshape mult}\,}
\newcommand{\Min}{\mbox{\upshape Min}\,}
\newcommand{\supp}{\mbox{\upshape supp}\,}

\newcommand{\arithdeg}{\mbox{\upshape arithdeg}\,}
\newcommand{\geomdeg}{\mbox{\upshape geomdeg}\,}

\newcommand{\MON}{\mbox{\upshape \MON}}
\newcommand{\Init}{\mbox{\upshape In}\,}
\newcommand{\sqfree}{\Phi}
\newcommand{\Star}{\mbox{\upshape st}\,}

\title{Symmetric iterated Betti numbers}

\author{Eric Babson}
\author{Isabella Novik}
\author{Rekha Thomas}
\thanks{Research partially supported by NSF grants DMS 0070571 and DMS
0100141} 
\address{Department of Mathematics, University of
  Washington, Seattle, WA 98195-4350, email:
  [babson,novik,thomas]@math.washington.edu}    
\date{\today}

\begin{document}

\begin{abstract}
  We define a set of invariants of a homogeneous ideal $I$ in a
  polynomial ring called the symmetric iterated Betti numbers of $I$.
  For $I_{\Gamma}$, the Stanley-Reisner ideal of a simplicial complex
  $\Gamma$, these numbers are the symmetric counterparts of the
  exterior iterated Betti numbers of $\Gamma$ introduced by Duval and
  Rose.  We show that the symmetric iterated Betti numbers of an ideal
  $I$ coincide with those of a particular reverse lexicographic
  generic initial ideal $\Gin(I)$ of $I$, and interpret these
  invariants in terms of the associated primes and standard pairs of
  $\Gin(I)$.  We verify that  for an ideal $I=I_\Gamma$
  the extremal Betti numbers of $I_\Gamma$ are 
  precisely the extremal  (symmetric or exterior) iterated Betti numbers of
  $\Gamma$.  We close with some results and conjectures about the
  relationship between symmetric and exterior iterated Betti numbers
  of a simplicial complex.
\end{abstract}
\maketitle

\section{Introduction}
The goal of this paper is to define and study a set of invariants of a
homogeneous ideal in a polynomial ring, called the {\em symmetric
iterated Betti numbers} of the ideal. For a simplicial complex
$\Gamma$, the symmetric iterated Betti numbers of the
Stanley-Reisner ideal of $\Gamma$ (also referred
to as the symmetric iterated Betti numbers of $\Gamma$) are preserved
by {\em symmetric algebraic shifting}.

We discuss two versions of algebraic shifting (both introduced by Kalai
\cite{BK}, \cite{K91}) which given a simplicial 
complex $\Gamma$ with vertex set $[n] :=
\{1,2,\ldots, n\}$ provide new simplicial complexes with the same vertex set.  
We denote these versions by 
$\Delta(\Gamma)$ for the {\em symmetric shifting}  of $\Gamma$
(see Definition \ref{shifting}) and by 
$\Delta^e(\Gamma)$ for the {\em exterior shifting}  of $\Gamma$
(see Definition \ref{eshifting}).  
For both of these operations it is known that: 
\begin{enumerate}
\item[(P1)] $\Delta^{(e)}(\Gamma)$ is {\em shifted}, that is, for every
  $F\in \Delta^{(e)}(\Gamma)$, if $j<i\in F$, then
  $(F\setminus \{i\})\cup \{j\} \in \Delta^{(e)}(\Gamma)$.
\item[(P2)]
   If $\Gamma$ is shifted, then $\Delta^{(e)}(\Gamma)=\Gamma$.
\item[(P3)]
    $\Gamma$ and $\Delta^{(e)}(\Gamma)$ have the same $f$-vector, that is,
they have the same number of $i$-dimensional faces for every $i$.
\item[(P4)]
 If $\Gamma'$ is a subcomplex of $\Gamma$, then
$\Delta^{(e)}(\Gamma')\subset \Delta^{(e)}(\Gamma)$.
\end{enumerate}

Both versions were studied extensively from the algebraic point of
view in a series of recent papers by Aramova, Herzog, Hibi and others
(surveyed in \cite{H}).

Consider the polynomial ring $\Sfield=\field[y_1, \ldots, y_n]$ where
$\field$ is a field of characteristic zero. Let $\mathbb N$ denote the
set of non-negative integers. If $A \subseteq [n]$ then write
$y^A=\prod_{a \in A}y_a$.  Denote by $\N^{[n]}$ the monomials of
$\Sfield$ by identifying a function $f:[n]\rightarrow \N$ in
$\N^{[n]}$ with the monomial $\prod_{i \in [n]}y^{f(i)}_{i}$ 
and consider $\N^{[n]}$ as a multiplicative monoid.  Thus
$\{0,1\}^{[n]}=y^{2^{[n]}}$ is the set of squarefree monomials.  If
$\Gamma \subseteq 2^{[n]}$ is a simplicial complex then the {\em
  Stanley-Reisner} ideal of $\Gamma$ \cite[Def.~II.1.1]{St} is the
squarefree monomial ideal
$$
I_{\Gamma} := \langle y^{2^{[n]}-\Gamma} \rangle \subset \Sfield.$$
The
(bi-graded) {\em Betti numbers} of a homogeneous ideal $I \subset
\Sfield$ are the invariants $\beta_{i,j}(I)$ that appear in the
minimal free resolution of $I$ as an $\Sfield$-module.  
\[ \rightarrow
\ldots \bigoplus_j \Sfield(-j)^{\beta_{i,j}(I)} \rightarrow \ldots
\rightarrow \bigoplus_j \Sfield(-j)^{\beta_{1,j}(I)} \rightarrow
\bigoplus_j \Sfield(-j)^{\beta_{0,j}(I)} \rightarrow I \rightarrow 0
\]
Here $\Sfield(-j)$ denotes $\Sfield$ with grading shifted by $j$.  We
say that $\beta_{i, i+j}(I)$, is {\em extremal} if $0
\neq\beta_{i,i+j}(I)=\sum_{i'\geq i, j'\geq j}\beta_{i',i'+j'}(I)$.
(This is equivalent to having $0 \neq\beta_{i,i+j}(I)$ and $0=
\beta_{i',i'+j'}(I)$ for every $i' \geq i$ and $j'\geq j$, 
$(i',j') \neq (i,j)$. )

Since $\Delta^{(e)}(\Gamma)$ are shifted complexes, their combinatorial
structures are simpler than that of $\Gamma$. Nonetheless,
$\Delta^{(e)}$ preserve many combinatorial and topological
properties.  
\begin{enumerate}
\item[1.]  $\Delta^{(e)}$ preserve topological
  Betti numbers: (see \cite[Thm.~3.1]{BK}, \cite[Prop.~8.3]{AH} for
  exterior shifting and \cite[Cor.~8.25]{H} for symmetric shifting).
  Moreover, exterior algebraic shifting preserves the {\em exterior
    iterated Betti numbers} of a simplicial complex. (There are two
  versions of exterior iterated Betti numbers --- one due to Kalai
  \cite[Cor.~3.4]{K93} and another due to Duval and Rose \cite{DR}.
  Both sets of numbers are preserved under exterior shifting.)
\item[2.]  $\Delta^{(e)}$ preserve Cohen-Macaulayness: a
  simplicial complex $\Gamma$ is Cohen-Macaulay if and only if
  $\Delta^{(e)}(\Gamma)$ is Cohen-Macaulay, which happens if and only if
  $\Delta^{(e)}(\Gamma)$ is pure (see \cite[Thm.~5.3]{K93},
  \cite[Prop.~8.4]{AH} for exterior shifting and \cite[Thm.~6.4]{K91}
  for symmetric shifting).
\item[3.]  $\Delta^{(e)}$ preserve extremal Betti numbers:
  $\beta_{i,i+j}(I_{\Gamma})$ is an extremal Betti number of
  $I_{\Gamma}$ if and only if $\beta_{i,i+j}(I_{\Delta^{(e)}(\Gamma)})$ is
  extremal for $I_{\Delta^{(e)}(\Gamma)}$, in which case
  $\beta_{i,i+j}(I_{\Gamma})=\beta_{i,i+j}(I_{\Delta^{(e)}(\Gamma)})$ (see
  \cite{BCP} for symmetric shifting and \cite[Thm.~9.7]{AH} for both
  versions.)
\end{enumerate} 
Property 3 is a far-reaching generalization of Property 2, while
Property 1 played a crucial role in Kalai's proof of Property 2 for
exterior shifting. This suggests that there might be a connection between
the iterated Betti numbers of a simplicial complex $\Gamma$ on the one
hand and the extremal Betti numbers of the ideal $I_{\Gamma}$ on the
other. This is one of the connections we establish in this paper.

Consider the action of $GL({\Sfield}_1)$ on $\Sfield$ and choose $u
\in GL({\Sfield}_1)$ to be generic.  Denote by $\m=\langle
{\Sfield}_1\rangle$ the irrelevant ideal of $\Sfield$.  If $I$ is a
homogeneous ideal in $\Sfield$ then write $J_0(I)=uI$ and
$J_i(I)=y_i\Sfield+(J_{i-1}(I): \m^\infty)$.  We now come to the
central definition of this paper.

\begin{definition}
The symmetric iterated Betti numbers of a homogeneous
 ideal $I$ in $\Sfield$ are 
\[
b_{i,r}(I):=\dim H^0(S/J_i(I))_r 
\quad \mbox{ for } 0\leq
i,r\leq n,
\]
where $H^0(-)_r$ stands for the $r$-th component of the 0-th local
cohomology with respect to the irrelevant ideal $\m$.  

If $\Gamma$ is a simplicial complex with vertex set $[n]$, define the
symmetric iterated Betti numbers of $\Gamma$ to be $b_{i,r}(\Gamma) :=
b_{i,r}(I_\Gamma)$, $0\leq i,r\leq n$.
\end{definition}

Our first result gives a combinatorial interpretation of the symmetric
iterated Betti numbers of a simplicial complex $\Gamma$ and shows that
they are invariant under symmetric algebraic shifting. Let
$\max(\Gamma)$ denote the set of facets (maximal faces) of $\Gamma$.
Write $\hbox{dim}(\Gamma)=\hbox{max}\{|F|-1 \, : \, F \in \Gamma\}$.

\begin{numtheorem}{\ref{local_cohomology}}
Let $\Gamma$ be a simplicial complex. Then 
\[
b_{i,r}(\Gamma)= 
  \left\{ \begin{array}{ll}
             |\{F\in \max(\Delta(\Gamma)) : |F|=i, \,
                        [i-r]\subseteq F,  \,\, i-r+1 \notin F\}|
                               &\mbox{ if }  r\leq i\\
                             0 &\mbox{ otherwise.}
        \end{array}
        \right.    
\]
In particular, since $\Delta(\Delta(\Gamma))=\Delta(\Gamma)$, it
follows that the symmetric iterated Betti numbers of $\Gamma$ are
invariant under symmetric shifting.
\end{numtheorem}

Theorem \ref{local_cohomology} implies that $b_{i,r}(\Gamma)=0$ unless
$0\leq r \leq i \leq \dim(\Gamma)+1$. The exterior iterated Betti
numbers of $\Gamma$, $b^e_{i,r}(\Gamma)$, defined by Duval and Rose
have precisely the same combinatorial formula (up to a slight change
in indices), except that in their definition, one replaces
$\Delta(\Gamma)$ by $\Delta^e(\Gamma)$ \cite[Thm.~4.1]{DR}.

The extremal Betti numbers of an ideal $I=I_\Gamma$
are the extremal  iterated Betti numbers (symmetric or exterior)
of the simplicial complex
$\Gamma$ in the following sense.

\begin{numtheorem}{\ref{Betti-extremal} and  \ref{ext-ext}}
Let $\Gamma$ be a simplicial complex.  The extremal Betti numbers of
$I_\Gamma$ form a subset of the symmetric as well as of the exterior
iterated Betti numbers of $\Gamma$. More precisely, $\beta_{j-1,
i+j}(I_\Gamma)$ is an extremal Betti number of $I_\Gamma$ if and only
if
$$b^{(e)}_{n-j',i'}(\Gamma)= 0 \quad \forall (i',j')\neq(i,j), \;
i'\geq i, \; j'\geq j, \mbox{ and } \, b^{(e)}_{n-j,i}(\Gamma)\neq
0.$$ In such a case $\beta_{j-1, i+j}(I_\Gamma)=b_{n-j,i}(\Gamma)=
b^{e}_{n-j,i}(\Gamma)$.
\end{numtheorem}

Let $\Gin(I)$ denote the reverse lexicographic {\em generic initial
  ideal} of a homogeneous ideal $I$ in $\Sfield$ with variables
  ordered as $y_n\succ y_{n-1} \succ \cdots \succ y_1$.  It follows
  from \cite[Cor.  1.7]{BCP} that the symmetric iterated Betti numbers
  of $I$ coincide with those of $\Gin(I)$.  We provide an alternate
  proof of this fact in Section \ref{section:local} (see Corollary
  \ref{BCPanalog}). Our next result interprets the symmetric iterated
  Betti numbers $b_{i,r}(I)$ in terms of the associated primes of
  $\Gin(I)$.  It is well known that all associated primes of $\Gin(I)$
  are of the form $P_{[i]} := \langle y_j \, : \, j \not \in [i]
  \rangle = \langle y_j \, : \, j > i \rangle$.

\begin{numtheorem}{\ref{primes}}
  The iterated Betti numbers of a homogeneous ideal $I$ are related to
  the ideal $\Gin(I)$. Those of an ideal $I_\Gamma$ are related to the
  ideals $\Gin(I_\Gamma)$, $I_{\Delta(\Gamma)}$, and the shifted
  complex $\Delta(\Gamma)$. The relationships are as follows.
\begin{enumerate}
\item The multiplicity of $P_{[i]}$ with respect to $\Gin(I)$ is
  $$\mult_{\Gin(I)}(P_{[i]}) = \sum_{r} b_{i,r}(I).$$
  If $I=I_\Gamma$ then 
  $$\mult_{\Gin(I_{\Gamma})}(P_{[i]}) = \sum_{r} b_{i,r}(\Gamma) = 
  | \{ F \in \max(\Delta(\Gamma)) \, : \, |F| =  i \}|.$$  
  
\item The degree, geometric degree, and arithmetic degree of
  $\Gin(I_{\Gamma})$ and $I_{\Delta(\Gamma)}$ have the following
  interpretations: 
\begin{eqnarray*}
&& (i) \qquad \deg(\Gin(I_{\Gamma})) = 
\geomdeg(\Gin(I_{\Gamma})) =
\sum_r {b_{d,r}(I_\Gamma)}
     \\
&&(i') \qquad = \deg(I_{\Delta(\Gamma)}) =  |\{ F
  \in \max(\Delta(\Gamma)) \, : \, |F| = d \}| \,\,; \nonumber \\
   \\
&&(ii) \qquad \arithdeg(\Gin(I_{\Gamma})) = 
   \sum_{i,r}b_{i,r}(I_\Gamma)  \\
&& (ii') \qquad = \arithdeg(I_{\Delta(\Gamma)}) = |\max(\Delta(\Gamma))|. 
\end{eqnarray*}
Equations (i) and (ii) also hold for arbitrary homogeneous
ideals $I$ in $\Sfield$.  
\end{enumerate}
\end{numtheorem}

This paper is organized as follows. In Section
\ref{section:preliminaries} we recall the basics of symmetric
shifting. Section \ref{section:monomial} defines and interprets
certain monomial sets that are at the root of all our proofs.  In
Sections \ref{section:local}--\ref{section:remarks} we prove the
theorems stated above. We conclude in Section \ref{section:remarks}
with some results and conjectures on the relationship between the
exterior and symmetric iterated Betti numbers of a simplicial complex.


\section{algebraic shifting}
\label{section:preliminaries} 
In this section we recall the basics of symmetric algebraic
shifting. (The description of exterior shifting is deferred to Section
7.)  For further details on symmetric and exterior shifting see the
survey articles by Herzog \cite{H} and Kalai~\cite{K00}.

Let $\N^\sigma$ denote the set of all finite degree monomials in the
variables $y_i$ with $i \in \sigma$ and $\N^\sigma_r$ denote the set
of elements of degree $r$ in $\N^{\sigma}$. In particular, if
$[n]=[1,n]=\{1, \ldots, n\}$ then $\N^{[n]}$ is the set of all
monomials in $S$ and $\{0,1\}^{\sigma}$ is the set of all square free
finite degree monomials in $\N^{\sigma}$.  In this paper we
fix the reverse lexicographic order $\succ$ on $\N^\Z$ with $y_{i}
\succ y_{i-1}$ for all $i \in \Z$ extending the partial ordering by
degree.  We also define the square free map $\Phi:\N^\Z \rightarrow
\{0,1\}^\Z$ to be the unique degree and order preserving bijection for
which $\Phi (y_0^n) = \prod_{-n<i\leq 0}y_i$.  Thus for example $\Phi
(y_4y_6^3y_7)=y_0y_3y_4y_5y_7$.

For each homogeneous ideal $I \subset S$ there exists a Zariski open
set $U(I) \subset GL(S_1)$ such that the ideal $\Init_{\succ}(u I)$,
(the {\em initial ideal} of $u I$ with respect to the monomial
order $\succ$ on $S$), is independent of the choice of $u \in U(I)$.
The ideal $\Init_{\succ}(u I)$ is called the {\em generic initial
  ideal} of $I$ with respect to $\succ$ and is denoted by
$\Gin(I)=\Gin_{\succ}(I)$ (see \cite[Chapter 15]{E}).  If $I$ is a
homogeneous ideal in $S$ then one way to explicitly and uniformly
construct an element $\alpha \in U(I)$ is to consider the extension
$\Field=\field(\{\alpha_{i,j}\}_{i,j \in [n]})/\field$ and then for
any ideal $I$ in $S$ the element
$$\alpha:S_\Field=S\otimes_\field \Field \rightarrow S_\Field \quad
\mbox{ given by } \quad \alpha y_i=\sum_{j=1}^n\alpha_{i,j}y_j$$
is
generic for $\Field I$ as an ideal of $S_\Field$.

For a homogeneous ideal $I$ in $S$ and a generic linear map $u \in
U(I)$ define $$B(I)= \{m \in \N^{[n]} \, : \, m \, \mbox{\small is not
in the linear span of } \{n|m\succ n\}\cup uI\}.$$ Note that $B(I)$ is
a basis of the vector space $M_0(I) = S/{u I}$ and hence
$B(I)=\N^{[n]}-\Gin(I)$.

\begin{definition}  \label{shifting}
  The symmetric algebraic shifting of a simplicial complex
  $\Gamma\subseteq 2^{[n]}$ is $\Delta(\Gamma)$ where 
$y^{\Delta(\Gamma)} =\Phi (B(I_\Gamma))\cap
  \N^{[1,\infty]}\subseteq \{0,1\}^{[n]}$.
\end{definition}

Note that this means that $I_{\Delta(\Gamma)}=\langle
\Phi(\N^{[n]}-B(I_\Gamma))\rangle$.

The fact that $\Delta(\Gamma)$ is a simplicial complex satisfying
conditions (P1)--(P4) was proved in \cite[Thm.~6.4]{K91}, \cite{AHH}
by using certain properties of $B(I)$.  We list some of them 
below.
\begin{enumerate}
\item[(B1)] $B(I)$ is a basis of $S/uI $, as well as of $S/\Gin(I)$.
\item[(B2)] $B(I)$ is an order ideal --- if $m\in
B(I)$ and $m'|m$, then $m'\in B(I)$.
\item[(B3)] $B(I)$ is shifted --- if $j<i$ and $y_im\in B(I)$ then 
$y_jm\in B(I)$.  
\end{enumerate}
(B1) was discussed above while (B2) follows from the fact that
$\Gin(I)$ is an ideal. (B3) is a consequence of the fact that generic
initial ideals are Borel fixed \cite[Theorem 15.20]{E}. In
characteristic 0, this is equivalent to $\Gin(I)$ being {\em strongly
  stable} \cite[Theorem 15.23]{E}, which means that if $j<i$ and 
$y_jm\in \Gin(I)$ then $y_im\in \Gin(I)$.  

In the case when $I=I_\Gamma$,
$B(I_\Gamma)$ has another fundamental property: 
\begin{enumerate}
\item[(B4)] If $m\in B(I_\Gamma)\cap \N^{[k,n]}_r$ and $r\geq k$ then 
$ m\N^{\{k\}}\subseteq B(I_{\Gamma})$ as well.
\end{enumerate}
This is due to Kalai \cite[Lemma 6.3]{K91}
and implies that $y_1, \ldots, y_n$ is an {\em almost regular}
$M_0(I_\Gamma)$-sequence (a notion
introduced by Aramova and Herzog \cite{AH}; it played a crucial
role in their proof that extremal Betti numbers are preserved by
algebraic shifting).

\section{Special monomial subsets} \label{section:monomial}
In this section we identify and interpret certain subsets of monomials
in the basis $B(I)$ of $M_0(I)$ that are at the root of
all our proofs. 

\begin{definition} \label{monomialsets}
Let $I$ be a homogeneous ideal in $\Sfield$.
  For $i \in [0,n]$ define
\begin{eqnarray*}
A_{i}(I)&:=& \left\{\begin{array}{lll} 
             m \in \N^{[i+1, n]} &: 
             & m\N^{\{i\}}\subseteq B(I),\,\,\,\, 
             m\N^{\{i+1\}} \not \subseteq B(I)\,\, 
             \end{array} \right\} \\
A_{i,r}(I) &:= & A_{i}(I)\cap\N^\N_r.
\end{eqnarray*}                                              
\end{definition}

Several remarks are in order. 
Since $B(I)$ is shifted (B3), $m\N^{\{i\}}\subseteq B(I)$ iff 
$m\N^{[i]}\subseteq B(I)$.  
Since $B(I_\Gamma)$ satisfies (B4), 
\begin{equation}  \label{r>i}
A_{i,r}(I_\Gamma)=\emptyset \qquad \mbox{ if } r>i \qquad \text{and hence} 
\qquad A_i(I_\Gamma) = \cup_{r=0}^{i} A_{i,r}(I_\Gamma).
\end{equation} 
Also if $m \in \N^{[r,n]}_r$ then 
$m\in B(I_\Gamma)$ iff $m\N^{[r]}\subseteq B(I_\Gamma)$.  
Hence
\begin{equation}  \label{A_i,r} 
A_{i,r}(I_\Gamma) = \left\{ 
             m \in \N_r^{[i+1, n]}  : 
                                      y_i^{i-r} \cdot m \in B(I_\Gamma),\
                  y_{i+1}^{i-r+1} \cdot m \not \in B(I_\Gamma) \right\}. 
\end{equation}

In \cite{STV}, Sturmfels, Trung and Vogel introduced a decomposition
of the standard monomials of an arbitrary monomial ideal $M$, called
its {\em standard pair decomposition}, in order to study the
multiplicities of associated primes and degrees of $M$. We study these
quantities for the monomial ideals $I_{\Gamma}$,
$I_{\Delta(\Gamma)}$, and $\Gin(I)$. In Section~\ref{section:primes}
we show their relationship to the symmetric iterated Betti numbers of
$\Gamma$ and $I$, respectively.  These results rely on the fact that
the sets of monomials $A_i(I)$ defined above index the {\em standard
pairs} of $\Gin(I)$. For a monomial $m\in\N^\Z$, let $\supp(m)
:= \{ i \, : \, y_i | m \}\subset \Z$ be called the {\em support} of $m$.
Thus $\supp:\{0,1\}^\sigma\rightarrow 2^\sigma$ is a bijection.  

\begin{definition} \label{standard_pairs} \cite{STV}
  Let $M=\langle M\cap \N^{[n]}\rangle \subseteq S$ be a monomial
  ideal.  A {\bf standard monomial} of $M$ is an element of
  $\N^{[n]}-M$.  An {\bf admissible pair} of $M$ is a subset
  $m\N^\sigma\subseteq \N^{[n]}-M$ with $m \in \N^{[n]-\sigma}$ or
  equivalently if we take $\Z^\sigma$ to be Laurent monomials then an
  admissible pair is a subset $m\Z^\sigma \cap \N^{[n]}$ with
  $m\Z^\sigma \cap M=\emptyset$.  A {\bf standard pair} of $M$ is a(n
  inclusion) maximal admissible pair.
\end{definition}

\begin{lemma}\label{stdpair_interpretation} 
  If $I\subseteq S$ is an ideal then the standard pairs of $\Gin(I)$
  are $\{a\N^{[i]}\ : \ a\in A_i(I)\}$.  (Here $[0]=\emptyset$.)
\end{lemma}
\proof We first argue that all standard pairs of $\Gin(I)$ are of the
form $a\N^{[i]}$ for some $i \in[0,n]$.  Suppose $m\N^\sigma$ is an
admissible pair of $\Gin(I)$ with $k=\hbox{max}(\sigma)$.  Since
$B(I)$ is shifted (B3) and $m\N^{\{k\}}\subseteq B(I)$ we obtain that
$m\N^{[k]}\subseteq B(I)$ and hence $m\N^\sigma\subseteq m\Z^{[k]}
\cap\N^{[n]}\subseteq B(I)$.  If $m\N^\sigma$ is standard
(maximal) this implies that $m\N^\sigma=m\Z^{[k]} \cap\N^{[n]}$ and
thus that $\sigma=[k]$.

If $m\N^{[i]}\subseteq B(I)$ is standard then by the above argument 
$m\N^{\{i+1\}}\not\subseteq B(I)$ so $m \in A_i(I)$.  

Finally, if $m \in A_i(I)$, then $m\N^{\{i\}}\subseteq B(I)$ 
and hence $m\N^{[i]}\subseteq B(I)$ is admissible.  
If $m\N^{[i]}\subseteq B(I)$ is not standard then 
$m\N^{[i]}\subset m'\N^{[i']}\subseteq B(I)$ so 
$i'>i$ and $m\N^{\{i+1\}}\subseteq B(I)$ contradicting the 
choice of $m$.  
\endproof

\begin{corollary}
If $m\N^{[i]}$ is a standard pair of $\Gin(I_\Gamma)$ then 
 the degree of $m$ is at most $i$. 
\end{corollary}

\begin{proof}
This follows from (\ref{r>i}).
\end{proof}

The standard pairs of monomial ideals of moderate size can be computed
using the computer algebra package {\em Macaulay~2} \cite{M2} (see the
chapter {\em Monomial Ideals} in \cite{EGSS} for details). This gives
a method for computing the sets $A_i(I)$ for small examples
--- see Example~\ref{ex} below. 

In the case when $I=I_\Gamma$ there is another
interpretation of the monomials in
$A_i(I_{\Gamma})$ that relates them to the shifted complex
$\Delta(\Gamma)$, and is useful for the proofs of Theorems
\ref{local_cohomology} and \ref{primes}. 

\begin{lemma} \label{Duvalbetti_interpretation}
There is a bijection between the sets
$$A_{i,r}(I_\Gamma)\; \mbox{ and }\; \{F\in \max(\Delta(\Gamma)) :
|F|=i, \, [i-r]\subseteq F, \, i-r+1 \notin F\}$$ given by $\sqfree$
with $A_{i,r}(I_\Gamma)\ni m \mapsto [i-r] \cup
\supp(\sqfree(m))=\supp(\sqfree(my_i^{i-r}))$.
\end{lemma}

\begin{proof}
  For $r>i$ the assertion follows from the fact that both sets are
  empty (see (\ref{r>i})).  To deal with the case $r\leq i$, note that
\begin{eqnarray}
&&\{\sqfree(m) \, :\,  m\in A_{i,r}(I_\Gamma)\} 
                   \stackrel{\mbox{\small by }(\ref{A_i,r})}{=}
\nonumber \\
&&  \left\{ \begin{array}{cc} \sqfree(m) \; : & 
                m \in B(I_\Gamma)\cap \N_r^{[i+1, n]}, \\
     & y_i^{i-r} \cdot m \in B(I_\Gamma), \,\, y_{i+1}^{i-r+1} \cdot
                            m \not \in B(I_\Gamma)
        \end{array}
                     \right\}\stackrel{\mbox{\small Def.}\ref{shifting}}{=} 
                                       \nonumber\\ 
&& \left\{ \begin{array}{cc}
         G \in \Delta(\Gamma) \, : & |G| = r, \, G \cap [i-r+1] = \emptyset, \\
               &   G \cup [i-r] \in \Delta(\Gamma), \,\, G \cup [i-r+1] \not
                                \in \Delta(\Gamma) 
                       \end{array}
               \right\}=        \label{set}\\
&& \left\{F\setminus [i-r] \, : \,
         F\in \max(\Delta(\Gamma)), \,\, |F|=i, \,\, 
            [i-r]\subseteq F, i-r+1 \notin F
       \right\}, \nonumber
\end{eqnarray} 
where in the last equality we used the fact that $\Delta(\Gamma)$ is
shifted.  Indeed, if $G\cup [i-r+1]\notin \Delta(\Gamma)$, then $G\cup
[i-r]\cup \{j\}\notin \Delta(\Gamma)$ for every $j>i-r+1$, $j\notin
G$, implying that $G\cup[i-r]$ is a facet of $\Delta(\Gamma)$ for
every element $G$ of the set (\ref{set}).
\end{proof}

\begin{corollary} \label{pairs-facets-bijection}
  The standard pairs of $\Gin(I_{\Gamma})$ are in bijection with the
  facets of $\Delta(\Gamma)$: $m\mathbb N^{[i]}$ is a standard pair of
  $\Gin(I_{\Gamma})$ if and only if $[i-r] \cup \supp(\sqfree(m))$ is
  a facet of $\Delta(\Gamma)$ of size $i$.
\end{corollary}

In Section~\ref{section:local} we verify that
$b_{i,r}(I)=|A_{i,r}(I)|$ for all $i, r \in [0,n]$, which via
Lemma~\ref{Duvalbetti_interpretation} proves Theorem
\ref{local_cohomology}. (Hence, in particular, it follows from Theorem
\ref{local_cohomology} that $A_i(I_{\Gamma}) = \emptyset$ for all $i >
\dim(\Gamma) + 1$.)  Thus the sets $A_{i,r}(I_\Gamma)$ and their
cardinalities $b_{i,r}(\Gamma)$ carry important information about
$\Gamma$, and we record them in the following {\em triangles}.

\begin{definition}
  The $b$-triangle and monomial $b$-triangle of a simplicial complex
  $\Gamma$ are the lower triangular matrices whose respective
  $(i,r)$-th entries are $b_{i,r}(\Gamma)$ and $A_{i,r}(I_\Gamma)$ for
  $0 \leq i \leq r \leq \dim(\Gamma)+1$.
\end{definition}

\begin{example} \label{ex}
Let $\Gamma$ be the simplicial complex whose facets are
\begin{eqnarray*}
\max(\Gamma)=\{ \{1,2,4\},\{1,2,6\},\{1,3,4\},\{1,3,7\},
\{1,5,6\},\{1,5,7\},\{2,3,5\}, \\ 
\{2,3,7\},\{2,4,5\},
\{2,6,7\},\{3,4,6\},\{3,5,6\},\{4,5,7\},\{4,6,7\}\}.
\end{eqnarray*}

\begin{figure}
\epsfig{file=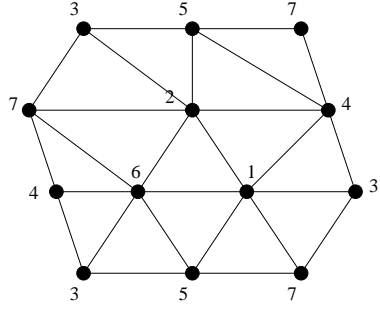, height=4cm}
\caption{The simplicial complex $\Gamma$ in Example~\ref{ex}.  
Here parallel boundary regions are identified.}
\label{f1}
\end{figure}
 
Then the Stanley-Reisner ideal of $\Gamma$ in the ring $S :=
 \field[a,b,c,d,e,f,g]$ is:
\begin{eqnarray*}
I_{\Gamma} = \langle
efg,cfg,afg,ceg,beg,cdg,bdg,adg,abg,def,bef,bdf, 
adf,bcf, \\
acf,cde,ade, 
ace, abe,bcd,abc \rangle. \end{eqnarray*}
Under the reverse lexicographic
order $\succ$ with $g \succ f \succ \cdots \succ b \succ a$, 
\begin{eqnarray*}
\Gin(I_\Gamma) &=& \langle
gf^2,f^3,f^2e,g^2f,gfe,fe^2,gfd,f^2d,fed,g^2e,ge^2,e^3,  \\
&& ged,e^2d,  fd^2,g^3,g^2d,
 gd^2,ed^2,g^2c,gfc,d^4 \rangle.
\end{eqnarray*}
Applying the map $\sqfree$  to the generators
  of $\Gin(I_{\Gamma})$ we get
\begin{eqnarray*}
I_{\Delta(\Gamma)} &= &\langle
gea,gfa,ecb,fcb,gcb,edb,fdb,gdb,feb,geb, 
gfb,edc,\\
&& fdc,gdc,fec,gec,gfc, fed,ged,gfd,gfe,dcba \rangle , 
\end{eqnarray*}
which shows that the shifted complex $\Delta(\Gamma)$ has facets:
\begin{eqnarray*}
\max(\Delta(\Gamma))   
&=& \{\{1,2,3\},\{1,2,4\},\{1,2,5\},\{1,2,
6\},\{1,2,7\}, \\
&& \{1,3,4\},
 \{1,3,5\}, 
\{1,3,6\},\{1,3,7\},
\{1,4,5\},\\
&&\{1,4,6\},\{1,4,7\},
\{1,5,6\},\{2,3,4\},\{5,
7\},\{6,7\}\}.
\end{eqnarray*}
The one skeleton of $\Delta(\Gamma)$ is (like that of $\Gamma$) 
the complete graph on it's $7$ vertices.  
The triangles are obtained by
coning $1$ with all edges involving the vertices
$2,3,4,5,6$ and $7$ except for $\{5,7\}$ and $\{6,7\}$ and adding the 
triangle $\{2,3,4\}$.
Thus all the triangles and the edges 
$\{5,7\}$ and $\{6,7\}$ are facets.  

We compute the $b$-triangle and the monomial $b$-triangle of $\Gamma$
by first computing the standard pairs of $\Gin(I_{\Gamma})$ using
Macaulay~2.

\begin{center}
\begin{tabular}{||l||l||}
\hline
$b$-triangle of $\Gamma$ & monomial
$b$-triangle of $\Gamma$ \\ 
\hline
\hline
$\begin{array}{l|llll}
  & 0 & 1 & 2 & 3 \\  
\hline
0 & 0 &   &   & \\
1 & 0 & 0 &   & \\
2 & 0 & 0 & 2 & \\
3 & 1 & 4 & 8 & 1
\end{array}$ & 
$\begin{array}{l|llll}
  & 0 & 1 & 2 & 3 \\  
\hline
0 & \emptyset &   &  & \\
1 & \emptyset & \emptyset &  & \\
2 & \emptyset & \emptyset & \{g^2, gf\} & \\
3 & \{1\} & \{g,f,e,d\} & \{ge, gd, f^2, fe, fd, e^2, ed, d^2 \} &
\{d^3\} \\
\end{array}$\\
\hline
\hline
\end{tabular}
\end{center}

The standard pairs of $I_{\Gamma}$, $\Gin(I_{\Gamma})$ and
$I_{\Delta(\Gamma)}$ are shown in the following table. Columns 2,3, and
  4 illustrate Lemma~\ref{Duvalbetti_interpretation} and
  Corollary~\ref{pairs-facets-bijection}.

\begin{center}
\begin{tabular}{|l|l|l|l|}
\hline
\upshape {StdPairs($I_{\Gamma}$)} & \upshape {StdPairs($\Gin(I_\Gamma
)$)} & $\sqfree$($m$) & \upshape{StdPairs($I_{\Delta(\Gamma)}$)}\\
form: $\mathbb N^\sigma$  & form: $m \mathbb N^{[i]}$ & $\rightarrow$
\supp($\sqfree(m)$) & form: $\mathbb N^\sigma$\\
\hline
\hline
$\mathbb N^{\{4,6,7\}}$ & $\mathbb N^{\{1, 2, 3\}}$   & $1$  $\rightarrow$
$\emptyset$ &  $\mathbb N^{\{1, 2, 3\}}$ \\ 
 $\mathbb N^{\{2,6,7\}}$ & $g\mathbb N^{\{1, 2, 3\}}$   & $g$  $\rightarrow$
 \{7\} & $\mathbb N^{\{1, 2, 7\}}$ \\  
$\mathbb N^{\{4,5,7\}}$ & $ge\mathbb N^{\{1, 2, 3\}}$  & $gd$ $\rightarrow$
\{7,4\} & $\mathbb N^{\{1, 4, 7\}}$ \\    
$\mathbb N^{\{1,5,7\}}$ & $gd\mathbb N^{\{1, 2, 3\}}$  & $gc$ $\rightarrow$
\{7,3\} & $\mathbb N^{\{1, 3, 7\}}$ \\  
$\mathbb N^{\{2, 3, 7\}}$ & $f\mathbb N^{\{1, 2, 3\}}$   & $f$ $\rightarrow$
\{6\} & $\mathbb N^{\{1, 2, 6\}}$ \\  
$\mathbb N^{\{1, 3, 7\}}$ & $f^2\mathbb N^{\{1, 2, 3\}}$ & $fe$ $\rightarrow$
\{6,5\} & $\mathbb N^{\{1, 5, 6\}}$ \\  
$\mathbb N^{\{3,5,6\}}$ & $fe\mathbb N^{\{1, 2, 3\}}$  & $fd$ $\rightarrow$
\{6,4\} & $\mathbb N^{\{1, 4, 6\}}$ \\  
$\mathbb N^{\{1,5,6\}}$ & $fd\mathbb N^{\{1, 2, 3\}}$  & $fc$ $\rightarrow$
\{6,3\} & $\mathbb N^{\{1, 3, 6\}}$ \\  
$\mathbb N^{\{3,4,6\}}$ & $e\mathbb N^{\{1, 2, 3\}}$   & $e$ $\rightarrow$
\{5\} & $\mathbb N^{\{1, 2, 5\}}$ \\  
$\mathbb N^{\{1,2,6\}}$ & $e^2\mathbb N^{\{1, 2, 3\}}$ & $ed$ $\rightarrow$
\{5,4\} & $\mathbb N^{\{1, 4, 5\}}$ \\  
$\mathbb N^{\{2,4,5\}}$ & $ed\mathbb N^{\{1, 2, 3\}}$  & $ec$ $\rightarrow$
\{5,3\} & $\mathbb N^{\{1, 3, 5\}}$ \\  
$\mathbb N^{\{2,3,5\}}$ & $d\mathbb N^{\{1, 2, 3\}}$   & $d$  $\rightarrow$
\{4\} & $\mathbb N^{\{1, 2, 4\}}$ \\  
$\mathbb N^{\{1,3,4\}}$ & $d^2\mathbb N^{\{1, 2, 3\}}$ & $dc$ $\rightarrow$
\{4,3\} & $\mathbb N^{\{1, 3, 4\}}$ \\  
$\mathbb N^{\{1,2,4\}}$ & $d^3\mathbb N^{\{1, 2, 3\}}$ & $dcb$ $\rightarrow$
\{4,3,2\} & $\mathbb N^{\{2, 3, 4\}}$ \\  
\hline
                  & $g^2\mathbb N^{\{1, 2\}}$    & $gf$ $\rightarrow$ \{7,6\}&
                  $\mathbb N^{\{6, 7\}}$ \\ 
                  & $gf\mathbb N^{\{1, 2\}}$     & $ge$ $\rightarrow$
                  \{7,5\}& $\mathbb N^{\{5, 7\}}$\\ 
\hline
\hline
\end{tabular}
\end{center}

\vspace{.5cm}
\end{example}

\section{Local cohomology} \label{section:local}
In this section we prove Theorem \ref{local_cohomology}, which provides 
a simple combinatorial formula for the symmetric iterated Betti
numbers of a simplicial complex.

\begin{theorem}\label{local_cohomology} 
For a simplicial complex $\Gamma$
\[
b_{i,r}(\Gamma) = 
\left\{ \begin{array}{ll}
             |\{F\in \max(\Delta(\Gamma)) :  |F|=i,
                        [i-r]\subseteq F,   i-r+1 \notin F\}| &\mbox{ if }
 r\leq i\\
             0 &\mbox{ otherwise}
        \end{array}
        \right.  
\]
\end{theorem}

The symmetric iterated Betti numbers $b_{i,r}(\Gamma)$ were defined as
the dimensions of the vector spaces $H^0(M_i(I_\Gamma))_r$, where for
a homogeneous ideal $I$ in $S$ and a generic linear map $u \in U(I)$,
$M_0(I)=S/uI $ and $M_i(I)=M_{i-1}(I)/(y_iM_{i-1}(I)+
H^0(M_{i-1}(I)))$ for $1\leq i\leq n$.  Thus at step $i$ we ``peel
off'' the $i$-th variable. This is similar to the ``deconing'' of the
shifted complex $\Delta(\Gamma)$ used in the definition of the
exterior iterated Betti numbers of $\Gamma$ by Duval and Rose
\cite{DR}.

In view of Lemma \ref{Duvalbetti_interpretation}, it suffices to show
that $|A_{i,r}(I)|=\dim H^0(M_i(I))_r$ for all $i,r\geq 0$ in order to
prove Theorem~\ref{local_cohomology}. We establish this in
Lemma~\ref{local} below. However, we first digress briefly to derive
and illustrate certain facts needed in the proof of Lemma~\ref{local}.

Recall that if $M$ is an $S$-module, $N$ is a submodule and $I$ is an
ideal in $S$ then $(N:I^\infty)_M=\{m\in M|\hbox{ for some } r\in \N,
I^r m\subseteq N\}$ and if $I = \langle f \rangle$ it is typical to
write $(N: \langle f \rangle^\infty) = (N:f^\infty)$.  For an
$S$-module $M$, the 0-th local cohomology of $M$ with respect to the
irrelevant ideal $\m=S_+=\langle y_1, \ldots, y_n \rangle$ is defined
as
$$
H^0(M)=\{m\in M : \m^k \cdot m=0 \mbox{ for some } k\}=(0:\m^\infty)_M.
$$
In particular $H^0(M)$ is graded when $M$ is graded. Hence the
equivalent definition of $M_i(I)$ is $M_i(I)=S/J_i(I)$, where
$$J_0(I)=u I \quad \mbox{ and } \quad
J_i(I)=y_iS+(J_{i-1}(I): \m^\infty).$$

Fix an $i$ such that $1 \leq i < n$. Then for all $1 \leq k < i$,
$(J_{i-1}(I) : y_k^\infty) = S$ since $y_k \in J_{i-1}(I)$.  For a $j$
such that $1 \leq i < j \leq n$, consider the family of automorphisms
$g_{a} \in \hbox{GL}(S_1)$ such that $g_a(y_i)=ay_i + (1-a)y_j$,
$g_a(y_j)=(1-a)y_i + ay_j$ and $g_a(y_r)=y_r$ otherwise, parameterized
by all $a\in \field$.  By induction $g_a(J_{i-1}(I))=J_{i-1}(I)$. Thus
$g_a(J_{i-1}(I) : y_i^\infty)= (J_{i-1}(I) : (ay_i+(1-a)y_j)^\infty)$,
and so the two colon ideals are isomorphic.  
If $I\subseteq S$ is a fixed ideal and $f \in S$ varies then 
the ideal $(I:f^\infty)$ depends only on which associated primes of 
$I$ contain $f$.  Thus if for a family of $f$'s over $\field$ all 
the colon ideals $(I:f^\infty)$ are isomorphic they must in fact be equal.  
Hence $$\forall \, i < j,\hbox{ we have } \, \, (J_{i-1}(I):
y_{i}^\infty)=(J_{i-1}(I):y_j^\infty)=(J_{i-1}(I): \m^\infty).$$

Fix $I$ and write $B=B(I)$, $A_i=A_i(I)$, $J_i=J_i(I)$ and $M_i =
M_i(I)$. For the proof of Lemma~\ref{local} we introduce the sets  
$$C_i=(B-\cup_{j \leq i}A_j)\cap \N^{[i+1,n]}.$$

\begin{lemma} \label{cifacts}
The sets $A_i$ and $C_i$ have the following properties:
\begin{enumerate}
\item $C_{i-1}$ is the disjoint union $C_{i-1} = A_i \,
\dot{\cup} \, C_i \, \dot{\cup} \, y_iC_{i-1}$.
\item $C_i$ and $C_i \cup A_i$ are shifted order ideals in
$\N^{[i+1,n]}$.
\item $\field C_i \cap J_i = \{0\}$.
\item $B \subseteq \field C_i+(J_i : y_{i+1}^\infty)$.
\end{enumerate}
\end{lemma}

\begin{proof} 
\paragraph{(1)} If $m\mathbb N^{[j]}$ is a standard pair of $\Gin(I)$ with $j \leq
i$, then no monomial in $m \mathbb N^{[j]}$ lies in $C_i$.  Thus $C_i$
is the set of all monomials in $B \cap \mathbb N^{[i+1,n]}$ that lie
in standard pairs of the form $\ast\mathbb N^{[j]}$ where $j > i$. In
other words, $C_i = \{ m \in B \cap \mathbb N^{[i+1,n]} \, : \, 
m\N^{\{i+1\}} \subseteq B \}$. This implies
that $y_{i+1}C_i \subseteq C_i$ and 
\begin{eqnarray*}
C_{i-1} & = & \{ m \in B \cap \mathbb N^{[i,n]} \, : \, m\N^{\{i\}}
        \subseteq B \}\\ & = & \{ m \in B \cap \mathbb N^{[i+1,n]} \,
        : \, m\N^{\{i+1\}} \subseteq B\}\,\, \dot{\cup}\\ &&\{ m \in B
        \cap \mathbb N^{[i+1,n]}\, : \,m\N^{\{i+1\}} \not\subseteq
        B\}\,\,\dot{\cup} \,\,y_iC_{i-1}\\ & = & C_i \, \dot{\cup} \, A_i
        \, \dot{\cup} \,y_iC_{i-1}
\end{eqnarray*}
\paragraph{(2)} The definitions of $C_i$ and $A_i$ and the fact that
        $B$ is shifted 
imply via an induction that $C_i$ and $C_i \cup A_i$ are shifted order
ideals in $\N^{[i+1,n]}$.
\smallskip
\paragraph{(3)} We establish this fact by induction on $i$. Note that
$\field C_0 \cap J_0 = \{0\}$ since all elements of $C_0\subseteq B$
are standard monomials of $\Gin(I) = \Init_{\succ}(J_0)$. Assume
$\field C_{i-1} \cap J_{i-1} = \{0\}$, but there exists $0 \neq f \in
\field C_i \cap J_i$. Since $f = \sum \beta_m m \in \field C_i$, each
$m\in C_i \subset \mathbb N^{[i+1,n]}$, and so $m \not \in \langle
y_1, \ldots, y_i \rangle$. Therefore, $f \in J_i = y_iS + (J_{i-1} \,
: \, y_i^\infty)$ implies that $f \in (J_{i-1} \, : \, y_i^\infty)$
--- i.e., $fy_i^k \in J_{i-1}$ for some $k$. Since $\field C_{i-1}
\cap J_{i-1} = \{0\}$, $fy_i^k = \sum \beta_m m \cdot y_i^k \not \in
\field C_{i-1}$ and we infer that at least one of the monomials $m$ is
not in $C_{i-1}$ (since $y_i C_{i-1} \subseteq C_{i-1}$).  However,
$C_i \subseteq C_{i-1}$, and hence this $m$ is also not in $C_i$,
which is a contradiction. Thus $\field C_{i} \cap J_{i} = \{0\}$.
\smallskip
\paragraph{(4)} Since $B\cap\N^{[i+1,n]}=(C_i\cup \cup_{j\leq
i}A_i)\cap \N^{[i+1,n]}$ 
it suffices to show that if $j\leq i$ then $A_j \subseteq \field
C_i+(J_i : y_{i+1}^\infty)$.  For every $a \in \cup_{j\leq i} A_j$,
there exists some $t > 0$ such that $ay_{i+1}^t \in \Gin(I) =
\Init_{\succ}(J_0) \subseteq \Init_{\succ}(J_i)$ (since $J_0 \subseteq
J_i$). Therefore there exists $f \in J_i$ such that $f = ay_{i+1}^t -
\sum \alpha_m m$ where $ay_{i+1}^t$ is the leading term of $f$ with
respect to $\succ$, $\alpha_m \in \field$, and the monomials $m$ are
standard monomials of $\Init_{\succ}(J_i)$. Hence $m \in B \cap
\mathbb N^{[i+1,n]}$ (since all standard monomials of
$\Init_{\succ}(J_i)$ are in $B$ and $y_1, \ldots, y_i \in
\Init_{\succ}(J_i)$). Further, since $ay_{i+1}^t \succ m$ for each $m$
and since $ay_{i+1}^t, m \in \N^{[i+1,n]}$, it follows that $y_{i+1}^t
| m$. Thus $a = f/y_{i+1}^t + \sum \alpha_m (m/y_{i+1}^t)$. Finally,
since each $m/y_{i+1}^t \prec a$ in the above sum, we are done by
induction.
\end{proof}

\begin{example}\label{ex2}
Consider the ideal $I = \langle z^6-5z^4y^2, z^3yx^3-3xyz^5, y^2z^2
\rangle \subset k[x,y,z]$. Under the reverse lexicographic order
$\succ$ with $z \succ y \succ x$, 
$$ \Gin(I) = \langle z^4, y^3z^3, y^5z^2, xy^4z^2, x^3y^2z^3, x^5yz^3
\rangle.$$
The standard pairs of $\Gin(I)$ are: 
\begin{eqnarray*}
&&\bullet \,\,\mathbb N^{\{1,2\}}, z\mathbb N^{\{1,2\}},\\
&&\bullet \,\,z^2 \mathbb N^{\{1\}}, yz^2 \mathbb N^{\{1\}}, y^2z^2
\mathbb N^{\{1\}},
y^3z^2 \mathbb N^{\{1\}}, z^3\mathbb N^{\{1\}},\\ 
&&\bullet \,\,y^4z^2\mathbb N^{\emptyset}, yz^3 \mathbb N^{\emptyset},
y^2z^3\mathbb N^{\emptyset}, xy^2z^3\mathbb N^{\emptyset}, 
x^2y^2z^3\mathbb N^{\emptyset},
xyz^3\mathbb N^{\emptyset}, x^2yz^3\mathbb N^{\emptyset}, \\ 
&& \indent x^3yz^3\mathbb N^{\emptyset},x^4yz^3\mathbb N^{\emptyset}. 
\end{eqnarray*}
Figure~\ref{f2} shows the decomposition of the standard monomials of
$\Gin(I)$ given by its standard pairs. The generators of $\Gin(I)$ are
the labeled black dots and the standard pair $y^m \mathbb N^{\sigma}$
is depicted by the cone $m +\mathbb R_{\geq 0}^\sigma$. Therefore,
\begin{eqnarray*}
&&A_2 = \{1,z\}, A_1 = \{z^2, yz^2, y^2z^2,y^3z^2,z^3\}\\
&&A_0 = \{y^4z^2, yz^3, y^2z^3, xy^2z^3, x^2y^2z^3, xyz^3, x^2yz^3,
x^3yz^3, x^4yz^3 \}\\
&& C_2 = \emptyset\\
&& C_1 = \mathbb N^{\{2\}} \cup z \mathbb N^{\{2\}}\\ 
&& C_0 = \mathbb N^{[1,2]} \cup z \mathbb N^{[1,2]} \cup z^2 \mathbb
N^{\{1\}}, yz^2 \mathbb N^{\{1\}}, y^2z^2 \mathbb N^{\{1\}},
y^3z^2\mathbb N^{\{1\}} ,z^3\mathbb N^{\{1\}}
\end{eqnarray*}
Note that each $C_{i-1}$ is the disjoint union $C_{i-1} = A_i \cup C_i
\cup y_iC_{i-1}$: 
\begin{eqnarray*}
&&C_2 = (A_3 = \emptyset) \cup (C_3 = \emptyset) \cup 
(zC_2 = z\emptyset)\\
&&C_1 = (A_2 = \{1,z\}) \cup (C_2 = \emptyset) \cup (yC_1 = y\mathbb
N^{\{2\}} \cup yz\mathbb N^{\{2\}})\\
&&C_0 = (A_1 = \{z^2, yz^2, y^2z^2,y^3z^2,z^3\}) \cup (C_1 = \mathbb
N^{\{2\}} \cup z \mathbb N^{\{2\}}) \cup \\
&& (xC_0 = x\mathbb N^{[1,2]}
\cup xz \mathbb N^{[1,2]} \cup xz^2 \mathbb 
N^{\{1\}}, xyz^2 \mathbb N^{\{1\}}, xy^2z^2 \mathbb N^{\{1\}},
xy^3z^2\mathbb N^{\{1\}} ,xz^3\mathbb N^{\{1\}})
\end{eqnarray*}

\begin{figure}
\input{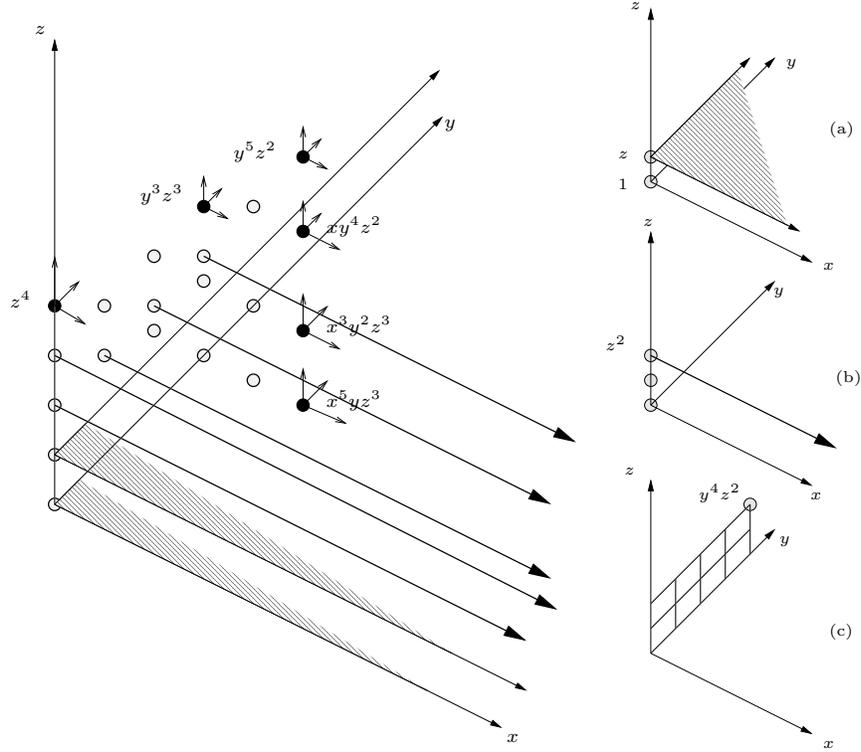}
\caption{The standard pair decomposition of $\Gin(I)$ where for
instance, the standard pairs $z\mathbb N^{\{1,2\}}$, $z^2\mathbb
N^{\{1\}}$ and $y^4z^2\mathbb N^{\emptyset}$ are shown separately 
in (a), (b) and (c) respectively.}
\label{f2}
\end{figure}
\end{example}

\begin{lemma}\label{local}
  $|A_{i,r}(I)|=\dim H^0(M_i(I))_r$ for all $i,r\geq 0$.
\end{lemma}

\begin{proof} Note that 
\begin{eqnarray*}
\dim H^0(M_i)_r & = & \dim(\{m \in (S/J_i)_r \, : \, m \in
                             (J_i \, : \, \m^\infty) \}) \\
                & = & \dim(\{m \in (S/J_i)_r \, : \, m \in
                             (J_i \, : \, y_{i+1}^\infty \}) \\
                & = & \dim (J_i : y_{i+1}^\infty)_r-\dim (J_i)_r.
\end{eqnarray*}

Hence, to prove the lemma it suffices to show that $\field A_i\oplus
J_i = (J_i:y_{i+1}^\infty)$. We do this by establishing the following
set of equalities:
$$
S \stackrel{[1]}{=}\field C_i\oplus \field A_i\oplus J_i
\stackrel{[2]}{=}\field C_i\oplus (J_i:y_{i+1}^\infty).
$$
$[1]_0$ follows from the facts that $\field B=\field A_0\oplus \field
C_0$ and $\field B\oplus J_0=S$.  We now show that $[1]_i$
implies $[2]_i$ and $[1]_{i+1}$.  Assume that $S=\field
C_i\oplus \field A_i\oplus J_i$.  Since $A_i \subseteq \field C_i+(J_i
: y_{i+1}^\infty)$ (Lemma~\ref{cifacts} (4)) we obtain that $S=\field
C_i+(J_i : y_{i+1}^\infty)$ and we must show that $\field C_i\cap (J_i
: y_{i+1}^\infty)=\{0\}$.  This follows from the facts that $y_{i+1}C_i
\subseteq C_i$ and that $\field C_i\cap J_i=\{0\}$
(Lemma~\ref{cifacts} (3)).  Thus
\begin{eqnarray*}
&&S \stackrel{[2]}{=}\field C_i\oplus (J_i : y_{i+1}^\infty)
\stackrel{\mbox{\small Lemma~\ref{cifacts} (1)}}{=}
\\
&&\field C_{i+1}\oplus\field A_{i+1}\oplus
   \field y_{i+1}C_{i}\oplus (J_i : y_{i+1}^\infty)
\stackrel{[1]}{=}\field C_{i+1}\oplus \field A_{i+1}\oplus J_{i+1}.
\end{eqnarray*}
\end{proof}

We close this section with several remarks.

\begin{remark}
Lemma~\ref{local} was verified in the special case of Stanley-Reisner
ideals of Buchsbaum complexes in \cite{N}.
\end{remark}

\begin{remark}
    Recall that by property (B1), 
    $B(I)$ is a basis of $S/uI$ as well as of $S/\Gin(I)$.
    Thus the proof of Lemma \ref{local} implies also that
    $$|A_{i,r}(I)|=\dim H^0({M}_i(\Gin(I))_r \quad \mbox{ for all } \quad i,
    r \geq 0,$$ and we recover the following fact (originally
    due to Bayer, Charalambous, and Popescu \cite[Cor.~1.7]{BCP}).
\end{remark}
\begin{corollary}   \label{BCPanalog}
Modules $H^0(M_i(I))$ and $H^0(M_i(\Gin(I)))$ 
have the same Hilbert function (for $i=0, 1, \ldots n$).
In other words, the symmetric iterated Betti numbers of $I$ are identical
to those of $\Gin(I)$.
\end{corollary}

\begin{remark}
Similar to the proof of Theorem \ref{local_cohomology}, one can show
that
\begin{eqnarray*}
&& \dim H^0(
S/\langle u I_\Gamma, y_1, \ldots, y_{i}\rangle)_r = \\
   &&          |\{G\in \Delta(\Gamma) :  |G|=r, \, [i-r+1]\cap G=\emptyset, \,
                        [i-r+1]\cup G\notin \Delta(\Gamma)\}|. 
\end{eqnarray*}
Thus for the shifted complex $\Gamma$, dimensions of the modules
$H^0(S/\langle u I_\Gamma, y_1, \ldots, y_{i}\rangle)_r$
coincide with Kalai's (exterior) iterated Betti numbers of $\Gamma$
(see \cite[Section 3]{K93}).  For that reason we refer to the numbers
$\bar{b}_{i,r}(\Gamma) :=\dim H^0(S/\langle u I_\Gamma,
y_1, \ldots, y_{i}\rangle)_r$ as Kalai's symmetric iterated Betti
numbers.
\end{remark}
\begin{remark}
Another fact worth  mentioning is that 
 $b_{i,i}(\Gamma)=\bar{b}_{i,i}(\Gamma)$ are just 
reduced (topological) Betti numbers of $\Gamma$, that is,
$$
b_{i,i}(\Gamma)=\bar{b}_{i,i}(\Gamma)= \beta_{i-1}(\Gamma) \quad
\forall \, 0\leq i\leq \dim(\Gamma)+1, \mbox{ where }
\beta_{i-1}(\Gamma)=\dim\widetilde{H}_{i-1}(\Gamma, \field).
$$
This result is a consequence of Theorem \ref{local_cohomology}
 together with the fact \cite{BK} that  for a shifted complex $K$ 
\[
\beta_{i-1}(K)=|\{F\in max(K) \, : \,|F|=i, \;1\notin F\}|,
\]
and the fact that symmetric shifting preserves topological Betti
numbers \cite{H}.  
\end{remark}
\begin{remark} \label{Buchsbaum}
Finally we note that if $\Gamma$ is a Buchsbaum complex
(i.e., a pure simplicial complex all of whose vertices have 
Cohen-Macaulay links), then
$b_{i,r}(\Gamma)={i-1 \choose r-1}\beta_{r-1}(\Gamma)$ for every 
$0\leq r\leq i\leq\dim(\Gamma)$, where $\beta_{r-1}(\Gamma)$
are reduced (topological) Betti numbers of $\Gamma$.
This follows from Lemma \ref{local} and \cite[Lemma 4.1]{N}.
\end{remark}
\section{Extremal Betti numbers} \label{section:betti}
This section is devoted to the proof of Theorem \ref{Betti-extremal},
which relates the graded algebraic Betti numbers of $I_{\Gamma}$
to the symmetric iterated Betti numbers of $\Gamma$.  Every
homogeneous ideal $I \subset \Sfield$ admits a graded free
$\Sfield$-resolution of the form
\[ \rightarrow
\ldots \bigoplus_j \Sfield(-j)^{\beta_{i,j}} \rightarrow \ldots
\rightarrow \bigoplus_j \Sfield(-j)^{\beta_{1,j}} \rightarrow
\bigoplus_j \Sfield(-j)^{\beta_{0,j}} \rightarrow I \rightarrow 0,
\]
where $\Sfield(-j)$ denotes $\Sfield$ with grading shifted by $j$.
Moreover, there exists a unique (up to isomorphism) resolution in
which all the exponents $\beta_{i,j}$ are simultaneously minimized,
called the {\em minimal graded free $\Sfield$-resolution} of $I$.
The numbers $\beta_{i,j}$ appearing in this minimal free resolution of
$I$ are called the {\em graded Betti numbers} of $I$.  A Betti number
of $I$, $\beta_{i, i+j}(I)$, is {\em extremal} if
$\beta_{i,i+j}(I)\neq 0$, but $\beta_{i', i'+j'}(I)=0$ for all $i'\geq
i$, $j'\geq j$, $(i',j')\neq (i,j)$. This terminology comes from the 
Betti diagram of $I$ output by the program {\em Macaulay~2} in
which the Betti numbers are arranged in a rectangular array whose 
columns are indexed by $i$ and rows by $j$ and the $(i,j)$-th entry is
the Betti number $\beta_{i,i+j}$. Thus $\beta_{i,i+j}$ is an extremal
Betti number of $I$ if it lies in a south-east corner of the {\em
  Macaulay~2} Betti diagram of $I$.

Let $\Gamma$ be a simplicial complex on the vertex set $[n]$.  The
{\em Alexander dual } of $\Gamma$ is the simplicial complex
\[
\Gamma^*=\{F\subseteq [n] : [n]\setminus F \notin \Gamma\}.
\]

The next two results (both due to Bayer, Charalambous and Popescu
\cite{BCP}, see \cite{AH} also for the second theorem) provide
connections between the extremal Betti numbers of the Stanley-Reisner
ideals of $\Gamma$ and $\Gamma^*$, and the shifted complex
$\Delta(\Gamma)$.
\begin{theorem}\label{extremal-dual}
  Let $\Gamma$ be a simplicial complex and $\Gamma^*$ be its
  Alexander dual.  The Stanley-Reisner ideals $I_{\Gamma}$ and
  $I_{\Gamma^*}$ have the same extremal Betti numbers. More precisely,
  $\beta_{i, i+j}(I_{\Gamma^*})$ is extremal if and only if
  $\beta_{j-1, i+j}(I_{\Gamma})$ is extremal. Also, in such a case
  $\beta_{i, i+j}(I_{\Gamma^*})=\beta_{j-1, i+j}(I_{\Gamma})$.
\end{theorem}

\begin{theorem}   \label{extremal-shifting}
  Extremal Betti numbers are preserved by algebraic shifting: for a
  simplicial complex $\Gamma$, $\beta_{i, i+j}(I_{\Gamma})$ is
  extremal if and only if $\beta_{i, i+j}(I_{\Delta(\Gamma)})$ is
  extremal. Moreover, in such a case $\beta_{i,
    i+j}(I_{\Gamma})=\beta_{i, i+j}(I_{\Delta(\Gamma)})$.
\end{theorem}

We are now in a position to prove the main theorem of this section.
\begin{theorem} \label{Betti-extremal}
The extremal Betti numbers of $I_{\Gamma}$ are contained among the 
  symmetric iterated Betti numbers of $\Gamma$. They are precisely the
  extremal entries in the $b$-triangle of $\Gamma$:
  $\beta_{j-1, i+j}(I_{\Gamma})$ is an extremal Betti number of
  $I_{\Gamma}$ if and only if
  $$b_{n-j',i'}(\Gamma)= 0 \quad \forall (i',j')\neq(i,j), \; i'\geq
  i, \; j'\geq j, \mbox{ and } b_{n-j,i}(\Gamma)\neq 0.$$
  Moreover, in this case,
  $\beta_{j-1,i+j}(I_{\Gamma})=b_{n-j,i}(\Gamma)$. 
\end{theorem}

\noindent{\bf Example~\ref{ex} continued:} The minimal free resolution
and Betti diagram of $I_{\Gamma}$ (computed by {\em Macaulay 2}) are
given below. Note that the entries in the southeast corners of the
Betti diagram of $I_{\Gamma}$ (the extremal Betti numbers of
$I_{\Gamma}$) are precisely the entries in the north-east corners of
the $b$-triangle of $\Gamma$ from Section~\ref{section:monomial}.

$$ 0 \rightarrow \Sfield^2 \rightarrow \Sfield^{15} \rightarrow 
\Sfield^{42} \rightarrow \Sfield^{49} \rightarrow
\Sfield^{21} \rightarrow \Sfield \rightarrow 0$$

$$\begin{array}{ccccccc}
\text{total:} &  1 & 21 & 49 & 42 & 15 & 2\\
         0:   &  1 & .  & .  &  . & .  & .\\
         1:   &  . & .  & .  &  . & .  & .\\
         2:   &  . & 21 & 49 & 42 & 14 & 2\\
         3:   &  . & .  & .  & .  & 1  & .
\end{array}$$

The proof of Theorem~\ref{Betti-extremal} relies on the following
lemma, which is a consequence of \cite[Thm.~2.1(b)]{HH} (see also
\cite[Prop.~12]{HRW}) and \cite[Cor.~6.2]{D}.  For completeness we
provide a different self-contained proof.

\begin{lemma} \label{Betti_algebraic}
The symmetric iterated Betti numbers of $\Gamma$ are related to the 
graded Betti numbers of the Stanley-Reisner
ideal $I_{\Delta(\Gamma^*)}$ as follows:
\[
\beta_{i, i+j}(I_{\Delta(\Gamma^*)})=
               \sum_r {n-r-j \choose i}b_{n-j,n-r-j}(\Gamma).
\]
\end{lemma}
\begin{proof}
  Let $\Gamma$ be a simplicial complex and let $\Gamma^*$ be its
  Alexander dual.  Recall that the Stanley-Reisner ideal
  $I_{\Delta(\Gamma^*)}\subset \Sfield$ is a squarefree monomial ideal
  whose minimal generators correspond to minimal non-faces of
  $\Delta(\Gamma^*)$. Let $G$ be the set of minimal generators of
  $I_{\Delta(\Gamma^*)}$, let $G_j = G \cap \mathbb
  N^{[n]}_j$, and let $\min(g)=\min\{i : y_i|g\}$ for a monomial
  $g\in G$.
 
  Since $\Delta(\Gamma^*)$ is a shifted complex, it follows that the
  ideal $I_{\Delta(\Gamma^*)}$ is {\em squarefree strongly stable},
  which means that for a monomial $m\in I_{\Delta(\Gamma^*)}$, if
  $y_i|m$, $i<j\leq n$, and $y_j$ is not a divisor of $m$, then
  $my_j/y_i \in I_{\Delta(\Gamma^*)}$ as well.  Hence the graded Betti
  numbers of $I_{\Delta(\Gamma^*)}$ are given by the following formula
  \cite[Cor.  3.4]{H} (which is the analog of the Eliahou-Kervaire
  formula for strongly stable ideals):
\begin{eqnarray}
&\beta_{i,i+j}(I_{\Delta(\Gamma^*)})&= \sum_{g\in G_j} {n-\min(g)+1-j
\choose i} \nonumber\\ &&=\sum_r {n-r-j \choose i}|\{g\in G_j :
\min(g)=r+1\}|.  \label{ii+j}
\end{eqnarray}

It is well known and is easy to prove that
$\Delta(\Gamma^*)=\Delta(\Gamma)^*$.  Thus, $g\in G_j$ if and only if
$\sigma=\supp(g)$ is a minimal non-face of $\Delta(\Gamma^*)$ of size
$j$, which happens if and only if $[n]\setminus\sigma$ is a facet of
$\Delta(\Gamma)$ of size $n-j$.  Moreover, $r+1=\min\{i :
i\in\sigma\}$ if and only if $[r]\subseteq [n]\setminus\sigma$, but
$r+1 \notin [n]\setminus\sigma$. Hence
\begin{eqnarray}
&&|\{g\in G_j : \min(g)=r+1\}| =\nonumber\\
&& |\{F\in\max(\Delta(\Gamma)) : 
         |F|=n-j, \; [r]\subseteq F, \mbox{ but } r+1 \notin F\}|
   \stackrel{\mbox{\small Thm }\ref{local_cohomology}}{=}
          \label{iterated}\\
&&b_{n-j, n-r-j}(\Gamma).\nonumber
\end{eqnarray}
Substituting (\ref{iterated}) in (\ref{ii+j}) gives the result.
\end{proof}

{\medskip\noindent {\it Proof of Theorem \ref{Betti-extremal}: \ }}
The theorem is an easy consequence of Lemma \ref{Betti_algebraic}. 
Indeed, since
${n-r-j \choose i}$ is positive for $r\leq n-i-j$ and is zero
otherwise, it follows from the lemma that
$$\beta_{i', i'+j'}(I_{\Delta(\Gamma^*)})=0 \quad \mbox{ iff } \quad
b_{n-j', n-j'-r}(\Gamma)=0 \mbox{ for all } r\leq n-i'-j'.
$$ 
Thus, 
\begin{eqnarray*}
&\beta_{i, i+j}(I_{\Delta(\Gamma^*)})\neq 0 \mbox{ is extremal }\quad   
\Longleftrightarrow &\\
&\beta_{i', i'+j'}(I_{\Delta(\Gamma^*)})=0 \mbox{ for all }i'\geq i,
\;j'\geq j, \;(i,j)\neq (i',j') \quad \Longleftrightarrow &\\  
&b_{n-j',i'}(\Gamma)=0 \mbox{ for all } i'\geq i,\; j'\geq j,\;
(i,j)\neq (i',j').& 
\end{eqnarray*}

Moreover, if this is the case, then all except the first summand in 
$$\beta_{i, i+j}(I_{\Delta(\Gamma^*)})={i \choose i} b_{n-j, i}(\Gamma)+
               \sum_{r<n-i-j} {n-r-j \choose i} b_{n-j, n-r-j}(\Gamma)
               $$
               vanish, implying that $\beta_{i,
                 i+j}(I_{\Delta(\Gamma^*)})=b_{n-j, i}(\Gamma)$.  
The result then follows from Theorems~\ref{extremal-dual} and 
\ref{extremal-shifting}.  \endproof

We remark that one can use \cite[Cor .~1.2]{AH} and certain properties
of sets $A_i(I)$ to provide a different proof of the following more
general result. We omit the details.
\begin{theorem} Let $I$ be a homogeneous ideal. The extremal Betti
numbers of $I$ form a subset of the symmetric iterated Betti numbers
of $I$. More precisely, $\beta_{j-1,i+j}(I)$ is an extremal Betti
number of $I$ if and only if
$$b_{n-j',i'}(I)= 0 \quad \forall (i',j')\neq(i,j), \; i'\geq
  i, \; j'\geq j, \mbox{ and } b_{n-j,i}(I)\neq 0.$$
  Moreover, in this case,
  $\beta_{j-1,i+j}(I)=b_{n-j,i}(I)$. 
\end{theorem}

\section{Associated Primes and Standard Pairs}\label{section:primes} 
The associated primes of a homogeneous ideal $I\subset S$ with
a primary decomposition $I = Q_1 \cap Q_2 \cap \cdots \cap Q_t$ are
the prime ideals $P_i := \sqrt{Q_i}$, $i = 1, \ldots, t$, where
$\sqrt{Q_i}$ denotes the radical of $Q_i$. The set of associated
primes of $I$, customarily denoted as $\Ass(I)$, is independent of the
primary decomposition of $I$. The minimal elements of $\Ass(I)$ with
respect to inclusion are called the {\em minimal primes} of $I$. We
denote the set of minimal primes of $I$ as $\Min(I)$.  Recall that the
irreducible (isolated) components of $V(I)$, the variety of $I$ in
$\field^n$, are the varieties $V(P)$ for $P\in\Min(I)$. Let
$Z_i:=V(P_i)$ be the variety of $P_i$ in $\field^n$. The finite
invariant $\deg(Z_i)$, called the {\em degree} of $Z_i$, is the
cardinality of $Z_i\cap L$ for almost all linear subspaces $L$ of
dimension equal to the codimension of $Z_i$.

\begin{definition} \label{def-degrees} \cite{BM}, \cite{STV}   
\begin{enumerate} 
\item
If $P$ is a homogeneous prime ideal in $S$ then the
{\bf multiplicity} of $P$ (with respect to $I$), denoted as $\mult_I(P)$ is
  the length of the largest ideal of finite length in the ring
  $S_P/IS_P$.
\item The {\bf degree} of $I$, $\deg(I) := \sum_{\{\dim(Z_i) =
    \dim(I)\}} \mult_I(P_i) \deg(Z_i).$
\item The {\bf geometric degree} of $I$, $$\geomdeg(I) :=
  \sum_{\{P_i \in \Min(I)\}} \mult_I(P_i) \deg(Z_i).$$ 
\item The {\bf arithmetic degree} of $I$, $$\arithdeg(I) :=
  \sum_{\{P_i \in Ass(I)\}} \mult_I(P_i) \deg(Z_i).$$
\end{enumerate}
\end{definition}

The invariant $\mult_I(P) > 0$ if and only if $P \in \Ass(I)$.  Our
main goal in this section is to prove Theorem~\ref{primes}.  We first
specialize Definition~\ref{def-degrees} to monomial ideals. If $M$ is
a monomial ideal, then every associated prime of $M$ is of the form
$P_{\sigma} := \langle y_j \, : \, j \not \in \sigma \rangle$ for some
set $\sigma \subseteq [n]$.  Hence $V(P_{\sigma})$ is the
$|\sigma|$-dimensional linear subspace spanned by $\{e_j : j \in
\sigma \}$ and $\deg(V(P_{\sigma})) = 1$. The three degrees of $M$
from Definition~\ref{def-degrees} are therefore appropriate sums of
multiplicities of ideals in $\Ass(M)$ with respect to $M$.

For a monomial ideal $M$ the multiplicities of associated primes as
well as all the degrees referred to in Definition~\ref{def-degrees}
can be read off from the standard pairs of $M$ (see Definition
\ref{standard_pairs}) as shown in the following lemma. The statements
in this lemma are either stated or can be derived easily from the
results in \cite{STV}.

\begin{lemma} \label{stdpair_facts} Let $M$ be a monomial
ideal. Then,\qquad 
\begin{enumerate}
\item the set of standard pairs of $M$ is well defined,
\item $*\N^\sigma$ is a standard pair of $M$ if and only if
  $P_{\sigma} \in \Ass(M)$,
\item $\N^\sigma$ is a standard pair of $M$ if and only if
  $P_{\sigma} \in \Min(M)$,
\item the dimension of $M$ is the maximal size of a set $\sigma$ such
  that $*\N^\sigma$ is a standard pair of $M$,
\item if $P_{\sigma} \in \Ass(M)$, then $\mult_M(P_{\sigma})$ is the 
  number of standard pairs of $M$ of the form $*\N^\sigma$ and
\item \begin{enumerate} 
  \item $deg(M)$ is the number of standard pairs $*\N^\sigma$ of
    $M$ such that $|\sigma| = dim(M)$,
  \item $geomdeg(M)$ is the number of standard pairs $*\N^\sigma$
    of $M$ such that $\N^\sigma$ is a standard pair of $M$ and 
  \item $arithdeg(M)$ is the total number of standard pairs of $M$.
  \end{enumerate}
\end{enumerate}
\end{lemma}

Lemma~\ref{stdpair_interpretation} showed that $m\N^\sigma$ is a
standard pair of $\Gin(I)$ if and only if $\sigma = [i]$ and
$m \in A_i(I)$ for some $0\leq i \leq n$. Combining this fact with
Lemma~\ref{stdpair_facts} we obtain the following.

\begin{corollary} \label{assprimes} (see also \cite[Corollary
15.25]{E}) \qquad
\begin{enumerate}
\item[(i)] $P_{[d]}$, $d = dim(I)$, is the unique minimal prime of
$\Gin(I)$ (if $I=I_\Gamma$ then $d=\dim \Gamma+1$), and \item[(ii)]
all embedded primes of $\Gin(I)$ are of the form $P_{[k]}$ for some $k
< d$.
\end{enumerate}
\end{corollary}

Thus the submonoids in the standard pairs of $\Gin(I)$
are { initial} intervals of $[n]$ while the cosets can be
complicated. On the other hand, for the square free monomial ideals
$I_{\Gamma}$ and $I_{\Delta(\Gamma)}$, the cosets of the standard 
pairs are trivial and the submonoids determine the ideals 
(cf.  Example~\ref{ex}).

\begin{corollary}  \label{squarefree_stdpairs}
If $\Gamma$ is a simplicial complex then $I_\Gamma= \bigcap_{\sigma
  \in \max(\Gamma)} P_{\sigma}$ is the irredundant prime decomposition
  of $I_\Gamma$. In particular, $I_\Gamma$ has no embedded primes and
  its standard pairs are $\{\N^\sigma \, : \, \sigma\in
  \max(\Gamma)\}$.
\end{corollary}

By Corollary~\ref{pairs-facets-bijection}, $m\N^{[i]}$ is a standard
pair of $\Gin(I_{\Gamma})$ if and only if $[i-r] \cup \supp(
\sqfree(m) )$ is a facet of $\Delta(\Gamma)$ of size $i$.
Combining this fact with Corollary~\ref{squarefree_stdpairs} we get
the following bijection as well.

\begin{corollary} \label{stdpair-stdpair-bijection}
There is a bijection between the standard pairs of $\Gin(I_{\Gamma})$ 
and those of $I_{\Delta(\Gamma)}$ given by: 
  $m\N^{[i]}$ is a standard pair of $\Gin(I_{\Gamma})$ with $\deg(m) =
  r$ if and only if $\N^{[i-r] \cup \supp(\sqfree(m))}$ is a standard
  pair of $I_{\Delta(\Gamma)}$. 
\end{corollary}

Theorem~\ref{primes} is now a corollary of the results in
Sections~\ref{section:monomial} and \ref{section:local}, and those
stated thus far in this section.

\begin{theorem}{\label{primes}}
  The iterated Betti numbers of a homogeneous ideal $I$ are related to
  the ideal $\Gin(I)$. Those of an ideal $I_\Gamma$ are related to the
  ideals $\Gin(I_\Gamma)$, $I_{\Delta(\Gamma)}$, and the shifted
  complex $\Delta(\Gamma)$. The relationships are as follows.
\begin{enumerate}
\item The multiplicity of $P_{[i]}$ with respect to $\Gin(I)$ is
  $$\mult_{\Gin(I)}(P_{[i]}) = \sum_{r} b_{i,r}(I).$$
  If $I=I_\Gamma$ then 
  $$\mult_{\Gin(I_{\Gamma})}(P_{[i]}) = \sum_{r} b_{i,r}(\Gamma) = 
  | \{ F \in \max(\Delta(\Gamma)) \, : \, |F| =  i \}|.$$  
  
\item The degree, geometric degree, and arithmetic degree of
  $\Gin(I_{\Gamma})$ and $I_{\Delta(\Gamma)}$ have the following
  interpretations: 
\begin{eqnarray*}
&& (i) \qquad \deg(\Gin(I_{\Gamma})) = 
\geomdeg(\Gin(I_{\Gamma})) =
\sum_r {b_{d,r}(I_\Gamma)}
     \\
&&(i') \qquad = \deg(I_{\Delta(\Gamma)}) = |\{ F
  \in \max(\Delta(\Gamma)) \, : \, |F| = d \}| \,\,; \nonumber \\
   \\
&&(ii) \qquad \arithdeg(\Gin(I_{\Gamma})) = 
   \sum_{i,r}b_{i,r}(I_\Gamma)  \\
&& (ii') \qquad = \arithdeg(I_{\Delta(\Gamma)}) = |\max(\Delta(\Gamma))|. 
\end{eqnarray*}
Equations (i) and (ii) also hold for arbitrary homogeneous
ideals $I$ in $\Sfield$.  
\end{enumerate}
\end{theorem}

\begin{proof}
(1) By Lemmas~\ref{stdpair_facts} (5)  
\begin{eqnarray*}
\mult_{\Gin(I_{\Gamma})}(P_{[i]}) & = & 
| \{\text{standard pairs of} \, \Gin(I) \, \text{of the form}
\, *\N^{[i]} \}|\\
& = & | A_i(I)| \qquad (\text{by Lemma~\ref{stdpair_interpretation}})\\
& = & \sum_r b_{i,r}(I) \qquad (\text{by Theorem~\ref{local}}).
\end{eqnarray*}
In particular, $P_{[i]}$ is an associated prime of $\Gin(I)$ if and
  only if $b_{i,r}(I) > 0$ for some $r$.  For a simplicial complex
  $\Gamma$ on $[n]$, by Lemma~\ref{Duvalbetti_interpretation},
  $$\mult_{\Gin(I_{\Gamma})}(P_{[i]}) = | A_i(I_{\Gamma})| = \sum_{r}
b_{i,r}(\Gamma) = | \{ F \in \max(\Delta(\Gamma)) \, : \, |F| =
i \}|.$$ In particular, $P_{[i]}$ is an associated prime of
$\Gin(I_{\Gamma})$ if and only if $\Delta(\Gamma)$ has a facet of
size $i$.\\

\noindent (2) The same lemmas along with Definition~\ref{def-degrees}
yield these results.
\end{proof}

We now establish certain further facts about $\Ass(\Gin(I_{\Gamma}))$. 

\begin{definition} \label{chainproperty}
  For any ideal $I$ in a polynomial ring, its poset of associated
  primes $\Ass(I)$ has the {\bf chain property} if whenever $P \in
  \Ass(I)$ is an embedded prime, then there exists $Q \in \Ass(I)$
  such that $P \supset Q$ and $\dim(Q) = \dim(P) + 1$.
\end{definition}

\begin{corollary}
  The poset $\Ass(\Gin(I_{\Gamma}))$ possesses the chain property if
  and only if $\max(\Delta(\Gamma))$ has the property that whenever
  $\Delta(\Gamma)$ has a facet of size $k\leq\dim(\Gamma)$ 
then it also has a facet of size $k+1$.
\end{corollary}

\begin{proof} By Corollary~\ref{assprimes}, the poset
  $\Ass(\Gin(I_{\Gamma}))$ has the chain property if whenever $P_{[k]}
  \in \Ass(\Gin(I_{\Gamma}))$ for some $k \leq \dim \Gamma$ then
  $P_{[k+1]} \in \Ass(\Gin(I_{\Gamma}))$. By Theorem~\ref{primes} (1),
  this is equivalent to the condition that whenever $\Delta(\Gamma)$
  has a facet of size $k \leq \dim \Gamma $ then it
  also has a facet of size $k+1$.
\end{proof}

\begin{corollary}
  If $\Gamma$ is a Buchsbaum complex (cf. Remark~\ref{Buchsbaum}) then
  $\Ass(\Gin(I_{\Gamma}))$ has the chain property.
\end{corollary}

\begin{proof}
  It was shown in \cite{N} that if $\Gamma$ is a $(d-1)$-dimensional
  Buchsbaum complex then $b_{i,r}(\Gamma) = {{i-1} \choose {r-1}}
  \beta_{r-1}$ for $i<d$, where $\beta_{r-1}$ is the reduced
  (topological) Betti number of $\Gamma$. Hence for $i < d$,
  $\mult_{\Gin(I_{\Gamma})}(P_{[i]}) = \sum_r b_{i,r} = \sum_r {{i-1}
  \choose {r-1}} \beta_{r-1}$.  Therefore if some $\beta_k \neq 0$
  then $\mult_{\Gin(I_{\Gamma})}(P_{[i]}) > 0$ for all $i \geq k$,
  which implies that $\Ass(\Gin(I_{\Gamma}))$ has the chain property.
\end{proof}

\section{Iterated Betti numbers: exterior versus symmetric} 
\label{section:remarks}
                                       
We close the paper with several remarks and conjectures on connections
between symmetric iterated Betti numbers and exterior iterated Betti
numbers of a simplicial complex. The superscript $e$ is used to
denote exterior shifting.

We start with a brief description of exterior algebraic shifting
extracted from \cite{H}.  Let $E=\bigwedge (\field[y_1, \ldots,
y_n]_1)=\bigwedge S_1$ be the exterior algebra over the
$n$-dimensional vector space $\Sfield_1$. A {\em monomial} in $E$ is
an expression of the form $m=y_{i_1}\wedge y_{i_2}\wedge \cdots \wedge
y_{i_k}$, where $1\leq i_1 < i_2 < \ldots < i_k \leq n$; the set
$\{i_1, i_2, \ldots, i_k\}$ is called the support of $m$, and is
denoted by $\supp(m)$. The {\em exterior Stanley-Reisner ideal} of a
simplicial complex $\Gamma$ on $[n]$ is
\[
J_\Gamma := \langle m\in E\, : \, m \mbox{ is a monomial}, \,
\supp(m)\notin \Gamma \rangle.
\]

\begin{definition}\label{eshifting}
  The exterior algebraic shifting of $\Gamma$, $\Delta^e(\Gamma)$, is
  the simplicial complex defined by
  $J_{\Delta^e(\Gamma)}:=\Gin(J_\Gamma)$, where $\Gin(J_\Gamma)$ is
  the generic initial ideal of $J_\Gamma$ with respect to the reverse
  lexicographic order with $y_n\succ y_{n-1}\succ \cdots \succ
  y_1$.
\end{definition}

The exterior iterated Betti numbers of a simplicial complex $\Gamma$
were introduced by Duval and Rose \cite{DR}. They have the following
combinatorial description (up to a slight change in the indexing), which
we adopt as their definition.

\begin{definition}
The exterior iterated Betti numbers of a simplicial complex $\Gamma$ are
\[
b^e_{i,r}(\Gamma)=|\{F\in\max(\Delta^e(\Gamma)) \,: \, 
|F|=i, \, [i-r]\subseteq F, \, i-r+1\notin F\}|.
\]
\end{definition}

Since $\Delta^e(\Delta(\Gamma))=\Delta(\Gamma)$ and
$\Delta(\Delta^e(\Gamma))=\Delta^e(\Gamma)$, the above definition
and Theorem~\ref{local_cohomology} imply that
\begin{eqnarray*}
b_{i,r}(\Gamma)&=& b^e_{i,r}(\Delta(\Gamma)) \\
 b^e_{i,r}(\Gamma)&=& b_{i,r}(\Delta^e(\Gamma)).
\end{eqnarray*}
Hence we infer the following corollary from Lemma \ref{Betti_algebraic}.

\begin{corollary}                        \label{sum}
For a simplicial complex $\Gamma$  
\[
\beta_{i, i+j}(I_{\Delta^e(\Gamma^*)})=
               \sum_r {n-r-j \choose i}b^e_{n-j,n-r-j}(\Gamma).
\]
\end{corollary}

Since exterior shifting preserves extremal Betti numbers \cite{AH},
the same proof as in  Theorem \ref{Betti-extremal}  yields
\begin{theorem}   \label{ext-ext}
$\beta_{j-1, i+j}(I_{\Gamma})$ is an extremal Betti number of 
$I_{\Gamma}$
if and only if 
$$b^e_{n-j',i'}(\Gamma)= 0 \quad \forall (i',j')\neq(i,j), \; i'\geq
i, \; j'\geq j, \mbox{ and } b^e_{n-j,i}(\Gamma)\neq 0.$$
Moreover, if
this is the case, then
$$\beta_{j-1, i+j}(I_{\Gamma})= \beta_{j-1,
  i+j}(I_{\Delta^e(\Gamma)})= \beta_{j-1,
  i+j}(I_{\Delta(\Gamma)})=b^e_{n-j,i}(\Gamma)=b_{n-j,i}(\Gamma).$$
\end{theorem}

We remark that there are no known connections between general
(non-extremal) graded Betti numbers of $I_{\Delta(\Gamma)}$ and
$I_{\Delta^e(\Gamma)}$.  However it is conjectured (see
\cite[Conj.~8.9]{H}) that the following holds.

\begin{conjecture} \label{Herzog}
For every simplicial complex $\Gamma$,
$\beta_{i,j}(I_{\Delta(\Gamma)})\leq \beta_{i,j}(I_{\Delta^e(\Gamma)})$.
\end{conjecture} 

In analogy, we propose the following. 

\begin{conjecture} \label{we}
For every simplicial complex $\Gamma$,
$b_{i,j}(\Gamma)\leq b^e_{i,j}(\Gamma)$.
\end{conjecture} 

Note that since all the coefficients in the expression of graded Betti
numbers in terms of iterated Betti numbers (see Lemma
\ref{Betti_algebraic} and Corollary \ref{sum}) are non-negative,
Conjecture \ref{we} if true would imply Conjecture \ref{Herzog}.

Conjecture \ref{Herzog} was verified by Aramova, Herzog, and Hibi
\cite{AHH2} in the case when $\Gamma$ is the Alexander dual of a
sequentially Cohen-Macaulay complex (a notion introduced by Stanley
\cite[Def.~II.2.9]{St}): they showed that in such a case
$$\beta_{i,j}(I_{\Delta(\Gamma)})=\beta_{i,j}(I_{\Delta^e(\Gamma)})=
\beta_{i,j}(I_\Gamma).$$

We have the following related result.
\begin{proposition} \label{seqC-M}
  Conjecture \ref{we} holds for all sequentially Cohen-Macaulay
  complexes.  More precisely, if $\Gamma$ is sequentially
  Cohen-Macaulay, then
  $$b_{i,r}(\Gamma)= b^e_{i,r}(\Gamma)=h_{i,r}(\Gamma),$$
  where
  $(h_{i,r}(\Gamma))_{0\leq r\leq i\leq \dim(\Gamma)+1}$ is the
  $h$-triangle of $\Gamma$.
\end{proposition}
The notion of $f$- and $h$-triangles was introduced by Bj\"orner and
Wachs \cite{BW}. We recall the definition.  For a simplicial complex
$\Gamma$ set
\[
f_{i,j}(\Gamma):=|\{F\in \Gamma\,:\, |F|=j,\, \dim(\Star F)=
i-1\}|,
\]
where $\Star F$ denotes the star of $F$ in $\Gamma$.
The $h$-triangle of $\Gamma$, 
$(h_{i,j}(\Gamma))_{0\leq j\leq i\leq \dim(\Gamma)+1}$, is defined by
\[
h_{i,j}(\Gamma)=\sum_{s=0}^j (-1)^{j-s} {i-s \choose j-s} f_{i,s}(\Gamma).
\]

{\smallskip\noindent {\it Proof of Proposition \ref{seqC-M}: \ }}
Duval \cite[Thm.~5.1]{D} showed that if $\Gamma$ is a sequentially
Cohen-Macaulay simplicial complex, then
$h_{i,j}(\Gamma)=h_{i,j}(\Delta^e(\Gamma))$. In his proof he
relied only on the  properties (P3) and (P4) of the operator $\Delta^e$,
and the fact that $\Gamma$ is Cohen-Macaulay if and only if $\Delta^e(\Gamma)$
is pure.
Since operator $\Delta$ (symmetric shifting) possesses all these properties as well,
 it follows that for a sequentially Cohen-Macaulay complex $\Gamma$,
\begin{equation}         \label{h-equality}
h_{i,j}(\Gamma)=h_{i,j}(\Delta^e(\Gamma))=h_{i,j}(\Delta(\Gamma)).
\end{equation}

Another result due to Duval \cite[Cor.~6.2]{D} is that for a shifted
complex $K$, $b^e_{i,j}(K)=h_{i,j}(K)$.  Thus
\begin{eqnarray}
h_{i,j}(\Delta^e(\Gamma))&=&
                           b^e_{i,j}(\Delta^e(\Gamma))=b^e_{i,j}(\Gamma)
                           \label{e} \quad \mbox{ and }\\
                           h_{i,j}(\Delta(\Gamma))&=&
                           b^e_{i,j}(\Delta(\Gamma))=b_{i,j}(\Gamma). \label{s}
\end{eqnarray}
Equations (\ref{h-equality})--(\ref{s}) imply the proposition.
\endproof

We close the paper with one additional conjecture for the special case
of Buchsbaum complexes.

\begin{conjecture}
  If $\Gamma$ is a Buchsbaum complex, then
  $b_{i,j}(\Gamma)=b^e_{i,j}(\Gamma)$, and hence
  $\beta_{i,j}(I_{\Delta(\Gamma)})=\beta_{i,j}(I_{\Delta^e(\Gamma)})$.
\end{conjecture}

\end{document}